\documentclass[11.5pt, a4paper]{article}
\usepackage[pagewise]{lineno}
\usepackage{graphicx}
\usepackage{amssymb}
\usepackage{latexsym, bm}
\usepackage{multicol}
\usepackage{indentfirst}
\usepackage{amssymb,amsfonts}
\usepackage{amsmath}
\usepackage{cases}
\usepackage[noadjust]{cite}
\usepackage[colorlinks,citecolor=red]{hyperref}
\newtheorem{theorem}{Theorem}
\newtheorem{lemma}{Lemma}

\newtheorem{definition}{Definition}

\usepackage{amssymb}
\numberwithin{equation}{section}
\numberwithin{theorem}{section}
\numberwithin{remark}{section}
\numberwithin{definition}{section}
\numberwithin{lemma}{section}
\numberwithin{corollary}{section}
\numberwithin{proposition}{section}
\title{The local characterizations of the singularity formation for the MHD equations}
\author{Wenke Tan\footnote{tanwenkeybfq@163.com}\quad Fan Wu\footnote{wufan0319@yeah.net}\\
{\small Key Laboratory of Computing and Stochastic Mathematics (Ministry of Education),}\\
{\small School of Mathematics and Statistics, Hunan Normal University,}\\
{\small Changsha, Hunan 410081, China}\\
}

\date{}

\begin{document}
\maketitle
{\bf Abstract:}
This paper characterizes the possible blow-up of solutions for the 3D magneto-hydrodynamics (MHD for short) equations. We first establish some $\epsilon$-regularity criteria in $L^{q,\infty}$ spaces for suitable weak solutions, and then together with an embedding theorem from $L^{p,\infty}$ space into a Morrey type
space to characterize the local behaviors of solutions near a potential singular point. More precisely, we show that if $z_{0}=\left(t_{0}, x_{0}\right)$ is a singular point, then for any $r>0$ it holds that
$$
\limsup _{t \rightarrow t_{0}^{-}}\left(\left\|u(t, x)-u(t)_{x_{0}, r}\right\|_{L^{3, \infty}\left(B_{r}\left(x_{0}\right)\right)}+\left\|b(t, x)-b(t)_{x_{0}, r}\right\|_{L^{3, \infty}\left(B_{r}\left(x_{0}\right)\right)}\right)>\delta^{*};
$$
$$
	\limsup\limits _{t \rightarrow t_{0}^{-}}\left(t_{0}-t\right)^{\frac{1}{\mu}} r^{\frac{2}{\nu}-\frac{3}{p}}\|(u,b)(t)\|_{L^{p, \infty}\left(B_{r}\left(x_{0}\right)\right)}>\delta^{*} \text { for } \frac{1}{\mu}+\frac{1}{\nu}=\frac{1}{2},\,2 \leq \nu \leq \frac{2 p}{3},\,  3<p\leq\infty;$$
$$
	\limsup\limits _{t \rightarrow t_{0}^{-}}\left(t_{0}-t\right)^{\frac{1}{\mu}} r^{\frac{2}{\nu}-\frac{3}{p}+1}\|(\nabla u,\nabla b)(t)\|_{L^{p}\left(B_{r}\left(x_{0}\right)\right)}>\delta^{*} \text { for } \frac{1}{\mu}+\frac{1}{\nu}=\frac{1}{2},\, \nu \in\left\{\begin{array}{ll}
		{[2, \infty],} & p\geq 3 \\
		{[2, \frac{2p}{3-p}],} & \frac{3}{2}\leq p<3
	\end{array}\right.
$$
where $\delta^{*}$ is a positive constant independent on $\nu$ and $p$.

\medskip

\section{Introduction}
This paper is concerned with the local behaviors near a possible singularity of suitable
weak solutions to the incompressible MHD  equations
in three spatial dimensions with unit viscosity and zero external force
\begin{equation}\label{1.1}
\left\{
             \begin{array}{lr}
             \partial_t u+(u\cdot\nabla)u-\Delta u+\nabla P=(b\cdot\nabla)b,& \\
            \partial_t b+(u\cdot\nabla)b-\Delta b=(b\cdot\nabla)u,&\\
             \nabla\cdot u=\nabla\cdot b=0,&\\
             u(x,0)=u_{0}(x),b(x,0)=b_{0}(x),&
\end{array}
\right.
\end{equation}
where $u$, $P$ and $b$ are quantities corresponding to the velocity of the fluid, its pressure and the magnetic field, and $u_{0},b_{0}$ are the prescribed initial data satisfying the compatibility condition $\nabla\cdot u_{0}=\nabla\cdot b_{0}$ = 0 in the distributional sense. Physically, equations \eqref{1.1} governs the dynamics of the velocity and magnetic fields in electrically conducting fluids, such as plasmas, liquid metals, and salt water. For more details on the related background, we refer the reader to \cite{CH}.

While many mathematicians have made a lot of important contributions to theoretical research of the incompressible 3D MHD equations \eqref{1.1}, many questions are left open. It is well known that the global existence of weak solutions, local
existence, and uniqueness of smooth solutions to equations \eqref{1.1} were established in \cite{DL,TEMAM}.
In the absence of magnetic field (i.e., $b=0$),  \eqref{1.1} reduces to the incompressible Navier-Stokes (NS) equations. The problem of the global regularity of solutions to the NS equations in three and higher space dimensions is a fundamental
question in fluid dynamics and is still widely open.
In a seminal paper \cite{LERAY}, Leray shown that for $3<p \leq \infty$, there exists $c_{p}$ such that the conditions
$$
\|u(t)\|_{L^{p}}<\frac{c_{p}}{(T-t)^{\frac{p-3}{2 p}}}
$$
imply the regularity of weak solutions on $[0, T]$.
 It is also possible to guarantee that a solution must remain smooth as long as there is control
on some scale critical quantity. The Prodi-Serrin and Escauriaza-Seregin-Sverak  regularity criterion \cite{GP, SJ,ELSGA} states that if
\begin{equation}\label{1.2}
u\in L^{q}\left(0,T;L^{p}\left(\Omega\right)\right) \quad \text{with} \quad\frac{2}{q}+\frac{3}{p}\leq1 \quad \text{and}\quad 3\leq p\leq\infty,
\end{equation}
then the solution is smooth and can be
continued to a smooth solution for some time $T'>T$. 
Recently, the endpoint case $p = 3$, $q = +\infty$ was further improved by Tao \cite{TAO}, who proved a triply logarithmic lower bound on the rate of blowup of the $L^{3}$ norm. Somewhat more precisely, Tao showed that for an absolute constant $c>0$, if $T<+\infty$,
\begin{equation}\label{1.3}
\limsup _{t \rightarrow T} \frac{\|u(\cdot, t)\|_{L^{3}}}{\left(\log \log \log \frac{1}{T-t}\right)^{c}}=+\infty.
\end{equation}
In \cite{KKH}, Kim-Kozono proved the local boundedness of a weak solution $u$ under the assumption that $\|u\|_{L^{r, \infty}\left([0, T] ; L^{s, \infty}\left(\mathbb{R}^{3}\right)\right)}$ is sufficiently small for some $(r, s)$ with $\frac{2}{r}+\frac{3}{s}=1$ and $3 \leq s<\infty .$ The limiting case of the regularity criteria obtained by Kim and Kozono was proved by He and Wang \cite{HW}, i.e., any weak solution for the NS or MHD equations is regular under the assumption that $\|u\|_{L^{2, \infty}\left([0, T] ; L^{\infty}\left(\mathbb{R}^{3}\right)\right)}$ is sufficiently small. These results involving the NS equations were improved by Wang-Zhang \cite{WZ}, they showed that $\left\|u_{3}\right\|_{L^{r, \infty}\left([0, T] ; L^{s, \infty}\left(\mathbb{R}^{3}\right)\right)} \leq M$ and $\left\|u_{h}\right\|_{L^{r, \infty}\left([0, T] ; L^{s, \infty}\left(\mathbb{R}^{3}\right)\right)} \leq \epsilon_{M}$ with $\frac{2}{r}+\frac{3}{s}=1$ and
$3<s \leq \infty$ imply the regularity of the suitable weak solution, where $u_{h}=\left(u_{1}, u_{2}\right), u=\left(u_{h}, u_{3}\right)$ and $\epsilon_{M}$ is a small constant depending on $M$. More recently, Tan-Yin \cite{TAN} extended previous results, they proved that $u$ is regular in $[-1,0] \times \mathbb{R}^{3}$ provided that $u$ satisfies
$$
\left\|u_{3}\right\|_{{L^{2, \infty}}\left([-1,0] ; B M O\left(\mathbb{R}^{3}\right)\right)}=M<\infty
$$
and
$$\left\|u_{h}\right\|_{{L^{2, \infty}}\left([-1,0] ; B M O\left(\mathbb{R}^{3}\right)\right)} \leq\epsilon.$$
Another important step towards a better understanding of the Navier-Stokes equations
is the partial regularity theory. Scheffer introduced the idea of suitable weak
solutions and generalized energy inequality in a series of papers \cite{SC,SC1,SC2}. Scheffer's
results were further generalized by Caffarelli, Kohn and Nirenberg \cite{CKN}, where a criterion
of $\epsilon$-regularity theorem was established and the nullity of one-dimensional parabolic
Hausdorff measure of the singular set was provided. Lin \cite{LFH} gave a greatly simplified
new proof by an indirect argument. Later, Choe and Lewis \cite{CHJ} improved the parabolic Hausdorff dimension by a logarithmic factor.

For incompressible MHD equations, we would like to mention here the paper of He-Xin \cite{HX} that first extended the Serrin type criteria of NS equations to the MHD equations only in terms of velocity, these results merit attention, especially, which indicates that the velocity field plays a more dominant role than the magnetic field on the regularity of solutions. 
For conditional regularity or blow-up criteria of NS and MHD equations, please refer to \cite{CW,CCT,LD,KKK,HBV,ME,WF,ZZF} for a more comprehensive discussions. On the other hand, He and Xin \cite{HX1} introduced the definition of suitable
weak solutions and obtained partial regularity theorems of suitable weak solutions.
Kang and Lee \cite{KL} found regularity criteria under the smallness of some scaling invariant
quantities for velocity field without a smallness condition for the magnetic field. But they
assumed uniform bounds of some scaling invariant quantities for magnetic field. And then, Wang-Zhang \cite{WZ1} presented a new interior regularity criterion for suitable weak solutions of the MHD equations that the scaling invariant quantity of velocity is involved only.

Investigation of singular weak Leray-Hopf solutions was first performed by Leray in \cite{LERAY}. In particular,  he proved that if $[0, T)$ is the maximal existence interval of a smooth solution $u$ for the 3D Navier-Stokes equations, then for any $3<p \leq \infty$, there exists a constant $c_{p}$ depending only on $p$ such that
$$
\|u(\cdot, t)\|_{L^{p}\left(\mathbb{R}^{3}\right)} \geq \frac{c_{p}}{(T-t)^{\frac{p-3}{2 p}}}.
$$
In general, except for the case of $p=3$ (see \cite{ELSGA,KN}), if the solution $u$ satisfies
$$\|u(\cdot, t)\|_{L^{p}\left(\mathbb{R}^{3}\right)} \leq \frac{C}{(T-t)^{\frac{p-3}{2 p}}}$$ with some constant $C>0$, the regularity of $u$ at time $t=T$ is still open. For the axi-symmetric Navier-Stokes equations, Chen-Strain-Yau-Tsai \cite{CSY,CSY1} and Koh-Nadirashvili-Seregin-Sver\'ak \cite{KN} made some important progresses, and established that the weak solution does not develop type-I singularities on time $t=T$, i.e.,
$$
\|u(x, t)\|_{L^\infty(\mathbb{R}^3)} \leq \frac{C}{\sqrt{T-t}}.
$$
Behaviour of the $L^{3}$ norm is more subtle. In a breakthrough paper, Escauriaza, Seregin and Sver\'ak showed that if $\left(\bar{x}, T^{*}\right)$ is a singular point then
\begin{equation}\label{1.4}
	\lim \sup _{t \rightarrow T^{*}}\|u(\cdot, t)\|_{L^{3}\left(B_{\delta}(\bar{x})\right)}=\infty
\end{equation}
for any fixed $\delta>0$.
Later in \cite{SG}, Seregin improved \eqref{1.4}:
\begin{equation}\label{1.5}
\lim _{t \rightarrow T^{*}}\|u(\cdot, t)\|_{L^{3}\left(\mathbb{R}^{3}\right)}=\infty.
\end{equation}
Recently in \cite{AB}, Albritton and Barker refined \eqref{1.4} and \eqref{1.5} to show that if $\Omega$ is a bounded domain with $C^{2}$ boundary one has
\begin{equation}\label{1.6} \lim _{t \rightarrow T^{*}}\|u(\cdot, t)\|_{L^{3}\left(B_{\delta}(\bar{x}) \cap \Omega\right)}=\infty
\end{equation}
 for any fixed $\delta>0$.

We may now go back to the issue of the ``concentration of solutions''. Such phenomenon was investigated for other equations e.g. nonlinear Schr\"odinger in the wake of the pioneering work \cite{MT}, also see \cite{HK,Hr}. In \cite{LOW}, for the NS equations, an interesting concentration result is proven for a weak Leray-Hopf solution $u$ which first blows up at $T^{*}>0$. In particular, they showed that there exists $t_{n} \rightarrow T^{*}$ and $x_{n} \in \mathbb{R}^{3}$ such that \begin{equation}\label{1.7}\left\|u\left(\cdot, t_{n}\right)\right\|_{L^{m}\left(B_{\sqrt{C(m)\left(T^{*}-t_{n}\right)}}\left(x_{n}\right)\right)} \geq \frac{C(m)}{\left(T^{*}-t_{n}\right)^{\frac{1}{2}\left(1-\frac{3}{m}\right)}}, \quad 3 \leq m \leq \infty.
\end{equation}
Later, Maekawa-Miura-Prange in \cite{MMP} improved
\eqref{1.7} and showed that for every $t \in\left(0, T^{*}\right)$ (not just a sequence $t_{n} \rightarrow T^{*}$) there exists $x(t) \in \mathbb{R}^{3}$ such that
\begin{equation}\label{1.8}\|u(\cdot, t)\|_{L^{p}\left(B_{\sqrt{C(m)\left(T^{*}-t\right)}}(x(t))\right)} \geq \frac{C(m)}{\left(T^{*}-t\right)^{\frac{1}{2}\left(1-\frac{3}{m}\right)}} \quad 3 \leq m \leq \infty.
\end{equation}
However, in \eqref{1.7} and \eqref{1.8} no information is provided on $x_{n}$ and $x(t)$. It is natural to ask whether the concentration phenomenon occurs on balls $B(x, R)$ with $R=O\left(\sqrt{T^{*}-t}\right)$ and with $\left(x, T^{*}\right)$ being a singular point. Recently, Barker-Prange in \cite{BP} gave this problem an affirmative answer for the $L^{3}$ for Leray-Hopf solutions which first blow-up at time $T^{*}$ and which satisfy the Type I bound:
\begin{equation}\label{1.9}
\sup _{\bar{x} \in \mathbb{R}^{3}} \sup _{r \in\left(0, r_{0}\right) }\sup_{T^{*}-r^{2}<t<T^{*}}\left(  \int_{B_{r}(\bar{x})}|u(x, t)|^{2} d x\right)^{\frac{1}{2}} \leq M
\end{equation}
 for a fixed radius $r_{0} \in(0, \infty]$ and $M, T^{*} \in(0, \infty)$. Then it holds
 $$
 \|u(\cdot, t)\|_{L^{3}\left(B_{R}(x)\right)} \geq \gamma_{\text {univ }}, \quad R=O(\sqrt{T-t}).
 $$
 In conclusion, the existed works characterized the singularity formation accurately by the critical $L^3$ norm.  But, if we consider the local characterization of singularity by $L^3$ norm, there is an unpleasant problem. On the one hand, It is well-known that if $|u(T,x)|\leq\frac{c}{|x|}$ with small enough $c$ then $(T,0)$ can not be a singular point. One the other hand, it is clear that $||u(T,\cdot)||_{L^3(B_r)}=\infty$ for any $r>0$. This means that one can not exclude such point from the singular set by using $L^3$ norm. Noticing that $||u(T,\cdot)||_{L^{3,\infty}(B_r)}=(\frac{4\pi}{3})^\frac{1}{3}c$, it is more natural to characterize the singularity formation for the Navier-Stokes equations by $L^{3,\infty}$ norm.
Very recently, the first author in \cite{TAN1} established the concentration rate for the $L^{p, \infty}$ norm of $u$ with $3 \leq p \leq \infty$
for the Navier-Stokes equations. Precisely, it showed that if $z_{0}=\left(t_{0}, x_{0}\right)$ is a singular point, then for any $r>0$ it holds that
\begin{equation}\label{1.10}
	\limsup _{t \rightarrow t_{0}^{-}}\left\|u(t, x)-u(t)_{x_{0}, r}\right\|_{L^{3, \infty}\left(B_{r}\left(x_{0}\right)\right)}>\delta^{*} \end{equation}
or
\begin{equation}\label{1.11}	\limsup\limits _{t \rightarrow t_{0}^{-}}\left(t_{0}-t\right)^{\frac{1}{\mu}} r^{\frac{2}{\nu}-\frac{3}{p}}\|(u)(t)\|_{L^{p, \infty}\left(B_{r}\left(x_{0}\right)\right)}>\delta^{*} \text { for } \frac{1}{\mu}+\frac{1}{\nu}=\frac{1}{2},\,2 \leq \nu \leq \frac{2 p}{3},\,  3<p\leq\infty,
\end{equation}
where $\delta^{*}$ is a positive constant independent on $p$ and $\nu$.

Concerning the collective behavior of the MHD equations \eqref{1.1}, to the best of our knowledge, there is no available result in this topic. In this paper, we shall consider the local characterization of singularity by the $L^{p,\infty}$ norm of $(u, b)$ with $3 \leq p \leq \infty$ and the $L^p$ norm of $(\nabla u, \nabla b)$ with $\frac{3}{2} \leq p \leq \infty$.

Before stating our main results, we first introduce the definition of suitable weak solutions and the notations of some scaled dimensionless quantities.

\begin{definition}
We say that a triple $(u,b,P)$ is a suitable weak solution of \eqref{1.1} in $\mathbb{R}^3\times
(0,T) $ if the following conditions hold:\\
1. $u, b \in L^{\infty}\left(0, T ; L^{2}\left(\mathbb{R}^3\right)\right) \cap L^{2}\left(0, T ; H^{1}\left(\mathbb{R}^3\right)\right)$; \\
2. $\eqref{1.1}_{1, 2}$ holds in the sense of distributions, that is,
$$
\int_{0}^{t}\left(u, \partial_{t} \phi+\Delta\phi+(u \cdot \nabla)\phi\right) d s+\left(u_{0}, \phi(0)\right)=\int_{0}^{t}\left((b\cdot\nabla)\phi, b\right) d s,
$$
$$
\int_{0}^{t}\left(b, \partial_{t} \phi+\Delta\phi+(u \cdot \nabla)\phi\right) d s+\left(u_{0}, \phi(0)\right)=\int_{0}^{t}\left((b\cdot\nabla)\phi, u\right) d s,
$$
for all $\phi\in C^{\infty}_{c}\left([0, T ) \times \mathbb{R}^3\right) $ with $\nabla\cdot\phi=0$, where $(\cdot, \cdot)$ is the scalar product in $L^2(\Omega)$;\\
3.  $(u, b)$ satisfies the energy inequality, that is,
$$
\|u(t)\|_{L^{2}}^{2}+\|b(t)\|_{L^{2}}^{2}+2 \int_{0}^{t}\|\nabla u(s)\|_{L^{2}}^{2}+\|\nabla b(s)\|_{L^{2}}^{2} d s \leq\left\|u_{0}\right\|_{L^{2}}^{2}+\left\|b_{0}\right\|_{L^{2}}^{2}, \quad 0 \leq t \leq T.
$$\\

4. the pressure $P \in L_{l o c}^{\frac{3}{2}}\left(\mathbb{R}^{3} \times(0, T)\right)$ and the following local energy inequality holds
\begin{equation}
\begin{split}\label{1.12}
&\int_{\mathbb{R}^{3}} (|u(x, t)|^2 +|b(x, t)|^2)\phi d x+2\int^t_0\int_{\mathbb{R}^{3}}(|\nabla u(x, t)|^2 +|\nabla b(x, t)|^2)\phi dx dt\\
\leq &\int_{0}^{t}\int_{\mathbb{R}^{3}} \left[(|u|^2 +|b|^2)(\Delta\phi+\partial_t\phi)+u\cdot\nabla\phi(|u|^2 +|b|^2+2P)-
(b\cdot u)(b\cdot\nabla\phi)\right]dxdt,
\end{split}
\end{equation}
for any $\phi\in C^{\infty}_{c}\left((0, T ] \times \mathbb{R}^{3}\right)$.
\end{definition}
We say a point $z_{0}=\left(t_{0}, x_{0}\right) \in(0, T] \times \mathbb{R}^{3}$ is a regular point of solution $(u,b)$ to \eqref{1.1} if there exists a non-empty neighborhood $\mathcal{O}_{z_{0}} \subset(0, T] \times \mathbb{R}^{3}$ of $z_{0}$ such that $u, b \in L^{\infty}\left(\mathcal{O}_{z_{0}}\right)$. The complement of the set of regular points will be called the singular set.

For the sake of completeness, we
introduce some invariant quantities: Let $(u, b, P)$ be a solution of the MHD equations \eqref{1.1}, let
\begin{equation}\label{1.13}
u_\lambda(x, t)=\lambda u(\lambda x, \lambda^2 t),\, b_\lambda(x, t)=\lambda b(\lambda x, \lambda^2 t),\, P_\lambda(x, t)=\lambda^2 P(\lambda x, \lambda^2 t)
\end{equation}
for any $\lambda>0$, then the family $(u_\lambda, b_\lambda, P_\lambda)$ is also a solution of \eqref{1.1}.
Define
$$
A(u, b, r, z)=\sup\limits _{t-r^2\leq s<t}\frac{1}{r}\int_{B_{r}(x)\times\{s\}}(|u|^2+|b|^2) dx,
$$
$$
B(u,b,r, z)=\frac{1}{r}\iint_{Q_{r}(z)}(|\nabla u|^2+|\nabla b|^2 ) dxdt,
$$
$$
C(u,b,r, z)=\frac{1}{r^2}\iint_{Q_{r}(z)}(|u|^3+|b|^3) dxdt,
$$
$$
D(P,r, z)=\frac{1}{r^2}\iint_{Q_{r}(z)}|P|^{\frac{3}{2}} dxdt,
$$
where
$$
B_{r}(x_0)=\{x\in \mathbb{R}^3: |x-x_0|<r\}, \,\, B_r=B_r(0), \,\, B=B_1;
$$
$$
Q_{r}(x_0)=B_{r}(x_0)\times (t_0-r^2, t_0), \,\, Q_r=Q_r(0), \,\, Q=Q_1.
$$
For simplicity, we denote
$$
A(u,b,r)=A(u,b,r, 0), \, B(u,b,r)=B(u,b,r, 0),\, C(u,b,r)=C(u,b,r, 0),\, D(P,r)=D(P,r, 0).
$$
%
Next, we recall the following function spaces:

We use $L^{q}\left((0, T] ; L^{p}\left(\mathbb{R}^{3}\right)\right)$ to denote the space of measurable functions with the following norm
$$
\|f\|_{L^{q}\left([0, T] ; L^{p}\left(\mathbb{R}^{3}\right)\right)}=\left\{\begin{array}{l}
	\left(\int_{0}^{T}\left(\int_{\mathbb{R}^{3}}|f(t, x)|^{p} d x\right)^{\frac{q}{p}} d t\right)^{\frac{1}{q}}, 1 \leq q<\infty \\
	\operatorname{ess} \sup _{t \in(0, T]}\|f(t, \cdot)\|_{L^{p}\left(\mathbb{R}^{3}\right)}, q=\infty
\end{array}\right.
$$
The Lorentz space $L^{r, s}([0, T])$ is the space of measurable functions with the following
norm:
$$
\|f\|_{L^{r, s}([0, T])}=\left\{\begin{array}{l}
	\left(\int_{0}^{\infty} \sigma^{s-1}|\{x \in[0, T]:|f(x)|>\sigma\}|^{\frac{s}{r}} d \sigma\right)^{\frac{1}{s}}, 1 \leq s<\infty \\
	\sup _{\sigma>0} \sigma|\{x \in[0, T]:|f(x)|>\sigma\}|^{\frac{1}{r}}, s=\infty
\end{array}\right.
$$
Throughout this paper, $u_{x_{0}, \rho} \doteq \frac{1}{\left|B_{\rho}\right|} \int_{B_{\rho}\left(x_{0}\right)} u~ d x$ and $C$ denotes an absolute positive number which can change from line to line.

With weak solutions and the relevant function spaces defined, we can now state the main
results of this paper.
\begin{theorem}\label{th1}Let $\frac{1}{q}+\frac{1}{p}=\frac{1}{2}$ with $2 \leq p \leq \infty .$ Assume $(u,b, P)$ be a suitable weak solution to the MHD equations \eqref{1.1} on $Q_{1}\left(z_{0}\right)$, there exists a absolute constant $\delta>0$ such that if
	\begin{equation}
	\begin{aligned}
		\left\| \left(\sup _{\eta \leq 1}\left(\frac{1}{\eta} \int_{B_{\eta}\left(x_{0}\right)}\left|u(x, t)-u_{x_{0}, \eta}\right|^{p} d x\right)^{\frac{1}{p}},  \sup _{\eta \leq 1}\left(\frac{1}{\eta} \int_{B_{\eta}\left(x_{0}\right)}\left|b(x, t)-b_{x_{0}, \eta}\right|^{p} d x\right)^{\frac{1}{p}}\right) \right\|_{L^{q, \infty}\left[t_{0}-1, t_{0}\right]} \leq \delta 	
	\end{aligned}
	\end{equation}
or
	\begin{equation}
	\begin{aligned}
	\left\|\left(\sup _{\eta \leq 1}\left(\frac{1}{\eta} \int_{B_{\eta}\left(x_{0}\right)}|u(x, t)|^{p} d x\right)^{\frac{1}{p}}, \sup _{\eta \leq 1}\left(\frac{1}{\eta} \int_{B_{\eta}\left(x_{0}\right)}|b(x, t)|^{p} d x\right)^{\frac{1}{p}}\right)\right\|_{L^{q, \infty}\left[t_{0}-1, t_{0}\right]} \leq \delta,
	\end{aligned}
	\end{equation}
	then $z_{0}$ is a regular point.
\end{theorem}
Let us now state the characterization for singularity formation in MHD equations \eqref{1.1}.
\begin{theorem}\label{th2} Let $(u, b, P)$ be a suitable weak solution of equations \eqref{1.1} in $Q_{1}\left(z_{0}\right)$. Assume $z_{0}$ be a singular point, then for any given $r \in(0,1)$, it holds that
	
\begin{equation}
\limsup _{t \rightarrow t_{0}^{-}}\left(\left\|u(t, x)-u(t)_{x_{0}, r}\right\|_{L^{3, \infty}\left(B_{r}\left(x_{0}\right)\right)}+\left\|b(t, x)-b(t)_{x_{0}, r}\right\|_{L^{3, \infty}\left(B_{r}\left(x_{0}\right)\right)}\right)>\delta^{*};
\end{equation}
\begin{equation}
	\begin{aligned}
	\limsup _{t \rightarrow t_{0}^{-}}\left(t_{0}-t\right)^{\frac{1}{\mu}} r^{\frac{2}{\nu}-\frac{3}{p}}\|(u,b)(t)\|_{L^{p, \infty}\left(B_{r}\left(x_{0}\right)\right)}>\delta^{*}, \frac{1}{\mu}+\frac{1}{\nu}=\frac{1}{2},  2 \leq \nu \leq \frac{2p}{3}, 3<p\leq\infty;
		\end{aligned}
\end{equation}
\begin{equation}
	\begin{aligned}
	\limsup\limits _{t \rightarrow t_{0}^{-}}\left(t_{0}-t\right)^{\frac{1}{\mu}} r^{\frac{2}{\nu}-\frac{3}{p}+1}\|(\nabla u,\nabla b)(t)\|_{L^{p}\left(B_{r}\left(x_{0}\right)\right)}>\delta^{*} \text { for } \frac{1}{\mu}+\frac{1}{\nu}=\frac{1}{2},\, \nu \in\left\{\begin{array}{ll}
	{[2, \infty],} & p\geq 3 \\
	{[2, \frac{2p}{3-p}],} & \frac{3}{2}\leq p<3
\end{array}\right.
\end{aligned}
\end{equation}
where $\delta^{*}>0$ is independent on $\mu, \nu, p$ and $r$.
\end{theorem}

\section{Several auxiliary lemmas}
In this section, we shall give several auxiliary lemmas, which present some uniform estimates of invariant quantities and a crucial embedding theorem from the Lorentz space $L^{p,\infty}$ to a Morrey type space.
\begin{lemma}\label{le2.1}
	Let $(u, b, P)$ be a suitable weak solution of equations \eqref{1.1} on $Q_1(z_0)$ satisfying
	\begin{equation}
	\begin{aligned}\label{2.1}
		&\left\| \left(\sup _{\rho\leq 1}\left(\frac{1}{\rho} \int_{B_{\eta}\left(x_{0}\right)}\left|u(x, t)-u_{x_{0},\rho}\right|^{p} d x\right)^{\frac{1}{p}},  \sup _{\rho \leq 1}\left(\frac{1}{\rho} \int_{B_{\eta}\left(x_{0}\right)}\left|b(x, t)-b_{x_{0}, \rho}\right|^{p} d x\right)^{\frac{1}{p}}\right) \right\|_{L^{q, \infty}\left[t_{0}-1, t_{0}\right]}\\
		& =M < \infty	
	\end{aligned}
\end{equation}
or
\begin{equation}
	\begin{aligned}\label{2.2}
		\left\|\left(\sup _{\rho \leq 1}\left(\frac{1}{\rho} \int_{B_{\rho}\left(x_{0}\right)}|u(x, t)|^{p} d x\right)^{\frac{1}{p}}, \sup _{\eta \leq 1}\left(\frac{1}{\rho} \int_{B_{\rho}\left(x_{0}\right)}|b(x, t)|^{p} d x\right)^{\frac{1}{p}}\right)\right\|_{L^{q, \infty}\left[t_{0}-1, t_{0}\right]}  =M < \infty,
	\end{aligned}
	\end{equation}
where $p$ and $q$ satisfy $\frac{1}{q}+\frac{1}{p}=\frac{1}{2}$,	then for $0<r<\rho\leq 1$, we have following estimates
	\begin{equation}\label{2.3}
		C\left(u,b, r, z_0\right) \leq C \frac{r}{\rho} C\left(u,b, \rho, z_0\right)+C\left(\frac{\rho}{r}\right)^{2} B^{\frac{9-3 p}{6-p}}\left(u,b, \rho, z_0\right) M^{\frac{3 p}{6-p}}, \,\,2 \leq p<3;
	\end{equation}
	\begin{equation}\label{2.4}
		C\left(u,b, r, z_0\right) \leq C \frac{r}{\rho} C\left(u,b, \rho, z_0\right)+C\left(\frac{\rho}{r}\right) A^{\frac{p-3}{p-2}}\left(u,b, \rho, z_0\right) M^{\frac{p}{p-2}},\,\, 3 \leq p \leq 6;
	\end{equation}
	\begin{equation}\label{2.5}
		C\left(u,b, r, z_0\right) \leq C \frac{r}{\rho} C\left( u,b, \rho, z_0\right)+C\left(\frac{\rho}{r}\right)^{\frac{3}{2}} A^{\frac{3}{4}}\left(u,b, \rho, z_0\right) M^{\frac{3}{2}}, \,\,6<p \leq \infty.
	\end{equation}
\end{lemma}
\textbf{Proof.}  First, we consider the assumption \eqref{2.1}, and the proof of conclusions will be
 divided into three cases: $$\text{case}\,1:2\leq p<3;$$$$\text{case}\,2:3\leq p\leq6;$$$$\text{case}\,3:6<p\leq \infty.$$
We define  $$f_{p}(t)=\left(\sup _{\rho \leq 1} \frac{1}{\rho} \int_{B_{\rho}\left(x_{0}\right)}\left|u(t, x)-u_{x_{0}, \rho}\right|^{p} d x\right)^{\frac{1}{p}} $$ and $$g_{p}(t)=\left(\sup _{\rho \leq 1} \frac{1}{\rho} \int_{B_{\rho}\left(x_{0}\right)}\left|b(t, x)-b_{x_{0}, \rho}\right|^{p} d x\right)^{\frac{1}{p}} .$$ For almost every time $t \in\left(t_{0}-\rho^{2}, t_{0}\right]$, we have
\begin{equation}\label{2.6}
	\begin{split}
		\int_{B_{r}\left(x_{0}\right)}\left(|u|^3 +|b|^3\right)dx
		&\leq\int_{B_{r}\left(x_{0}\right)}\left(|u-u_{x_{0}, \rho}|^3+|b-b_{x_{0}, \rho}|^3\right) dx+C|B_{r}|\left(|u_{x_{0}, \rho}|^3+|b_{x_{0}, \rho}|^3\right)\\
		&=J_1+J_2.
	\end{split}
\end{equation}
For $\text{case}\,1:2\leq p<3$, $J_1$ can be estimated as
\begin{equation}\label{2.16}
	\begin{split}
		J_{1} &\leq C\left\|u-u_{x_0,\rho}\right\|_{L^{p}}^{\frac{3 p}{6-p}}\left\|u-u_{x_0, \rho}\right\|_{L^{6}}^{\frac{6(3-p)}{6-p}}+ C\left\|b-b_{x_0,\rho}\right\|_{L^{p}}^{\frac{3 p}{6-p}}\left\|b-b_{x_0, \rho}\right\|_{L^{6}}^{\frac{6(3-p)}{6-p}}\\
		&\leq C\left(\left\|u-u_{x_0, \rho}\right\|_{L^{6}}^{\frac{6(3-p)}{6-p}}+\left\|b-b_{x_0, \rho}\right\|_{L^{6}}^{\frac{6(3-p)}{6-p}}\right)\left(\left\|u-u_{x_0,\rho}\right\|_{L^{p}}^{\frac{3 p}{6-p}}+\left\|b-b_{x_0,\rho}\right\|_{L^{p}}^{\frac{3 p}{6-p}}\right).
	\end{split}
\end{equation}
Integrating with respect to time from $t_0-r^{2}$ to $t_0$ and using Holder's inequality, one has
\begin{equation}\label{2.8}
	\begin{aligned}
		& \int_{t_0-r^{2}}^{t_0} \int_{B_{r}\left(x_{0}\right)}J_1 d x d s \\
		\leq &C\left(\int_{t_0-r^{2}}^{t_0}\left\|u-u_{x_{0}, \rho}\right\|_{L^{6}\left(B_{\rho}\left(x_{0}\right)\right)}^{2}+\left\|b-b_{x_{0}, \rho}\right\|_{L^{6}\left(B_{\rho}\left(x_{0}\right)\right)}^{2} d s\right)^{\frac{9-3 p}{6-p}}\\
		&\times\left(\int_{t_0-r^{2}}^{t_0}\left\|u-u_{x_{0}, \rho}\right\|_{L^{p}\left(B_{\rho}\left(x_{0}\right)\right)}^{\frac{3 p}{2p-3}}+\left\|b-b_{x_{0}, \rho}\right\|_{L^{p}\left(B_{\rho}\left(x_{0}\right)\right)}^{\frac{3 p}{2p-3}} d s\right)^{\frac{2 p-3}{6-p}} \\
		\leq &C\left(\int_{t_0-r^{2}}^{t_0} \int_{B_{\rho}\left(x_{0}\right)}(|\nabla u|^{2}+|\nabla b|^{2}) d x d s\right)^{\frac{9-3 p}{6-p}} \rho^{\frac{3}{6-p}}\left(\int_{t_0-r^{2}}^{t_0}f_p^{\frac{3 p}{2 p-3}}+g_p^{\frac{3 p}{2 p-3}} d s\right)^{\frac{2 p-3}{6-p}}.
	\end{aligned}
\end{equation}
Using the assumption \eqref{2.1}, and if $2<p<3$, we see that
\begin{equation}\label{2.9}
	\begin{aligned}
		 \int_{t_0-r^{2}}^{t_0}f_p^{\frac{3 p}{2 p-3}}+g_p^{\frac{3 p}{2 p-3}} d s\leq& \int_{t_0-r^{2}}^{t_0}\left(f_p+g_p\right)^{\frac{3 p}{2 p-3}} d s\\
		  =& \frac{3 p}{2 p-3} \int_{0}^{\infty} \sigma^{\frac{3+p}{2 p-3}}\left|\left\{s \in\left[t-r^{2}, t\right] ;\left(f_p+g_p\right)>\sigma\right\}\right| d \sigma \\
		\leq & \frac{3 p}{2 p-3}\left\{\int_{0}^{R} \sigma^{\frac{3+p}{2 p-3}} r^{2} d \sigma+M^{\frac{2 p}{p-2}} \int_{R}^{\infty} \sigma^{\frac{3+p}{2 p-3}-\frac{2 p}{p-2}} d \sigma\right\} \\
		\leq & R^{\frac{3 p}{2 p-3}} r^{2}+\frac{3(p-2)} {p}{R}^{\frac{3 p}{2 p-3}-\frac{2 p}{p-2}} M^{\frac{2 p}{p-2}} \\
		\leq &\left(\frac{4p-6}{p}\right) r^{\frac{p}{2 p-3}} M^{\frac{3 p}{2 p-3}},
	\end{aligned}
\end{equation}
where we take $R=r^{-\frac{p-2}{p}} M$. When $p=2$, we take $R=M$ such that $$\mu{\left\{s \in\left[t-r^{2}, t\right] ;\left(f_p+g_p\right)>M\right\}}=0.$$
Indeed,
\begin{equation}\label{2.10}
	\begin{aligned}
		\int_{t_0-r^{2}}^{t_0}f_p^{\frac{3 p}{2 p-3}}+g_p^{\frac{3 p}{2 p-3}} d s\leq& \int_{t_0-r^{2}}^{t_0}\left(f_p+g_p\right)^{\frac{3 p}{2 p-3}} d s\\
		=& \frac{3 p}{2 p-3} \int_{0}^{\infty} \sigma^{\frac{3+p}{2 p-3}}\left|\left\{s \in\left[t-r^{2}, t\right] ;\left(f_p+g_p\right)>\sigma\right\}\right| d \sigma \\
		\leq & \frac{3 p}{2 p-3}\int_{0}^{M} \sigma^{\frac{3+p}{2 p-3}} r^{2} d \sigma
		\leq \left(\frac{4p-6}{p}\right) r^{\frac{p}{2 p-3}} M^{\frac{3 p}{2 p-3}},
	\end{aligned}
\end{equation}
where we used the fact $r<\rho\leq 1$.
The estimate of \eqref{2.9}-\eqref{2.10} together with \eqref{2.8} yield that
\begin{equation}\label{2.11}
	\begin{aligned}
		& \int_{t_0-r^{2}}^{t_0} \int_{B_{r}\left(x_0\right)}J_1 d x d s \\
		\leq &C\left(\frac{4p-6}{p}\right)^{\frac{2 p-3}{6-p}} \rho^{\frac{3}{6-p}} r^{\frac{p}{6-p}}\left(\int_{t_0-r^{2}}^{t_0} \int_{B_{\rho}\left(x_0\right)}(|\nabla u|^{2}+|\nabla b|^{2}) d x d s\right)^{\frac{9-3 p}{6-p}}M^{\frac{3 p}{6- p}}.
	\end{aligned}
\end{equation}
For $J_2$, we have
\begin{equation}\label{2.12}
	\begin{aligned}
		J_2\leq C\left(\frac{r}{\rho}\right)^{3} \int_{B_{\rho}(x_0)}\left(|u|^{3}+|b|^{3}\right) d x
	\end{aligned}
\end{equation}
Combining \eqref{2.6}, \eqref{2.11} and \eqref{2.12}, one get
\begin{equation}\label{2.13}
	\begin{split}
		\frac{1}{r^2}\int_{Q_{r}(z_0)}\left(|u|^3 +|b|^3\right)dx dt
		\leq &C \rho^{\frac{3}{6-p}} r^{\frac{3p-12}{6-p}}\left(\int_{t-r^{2}}^{t} \int_{Q_{\rho}\left(z_0\right)}(|\nabla u|^{2}+|\nabla b|^{2}) d x d s\right)^{\frac{9-3 p}{6-p}}M^{\frac{3 p}{6- p}}\\
		&+C\frac{r}{\rho}
		\frac{1}{\rho^2} \int_{\mathbb{P}_{\rho}(z)}\left(|u|^{3}+|b|^{3}\right) d x dt.
	\end{split}
\end{equation}
That is to say
\begin{equation}\label{2.14}
	C\left(u,b, r, z_0\right) \leq C \frac{r}{\rho} C\left(u,b, \rho, z_0\right)+C\left(\frac{\rho}{r}\right)^{2} B\left(u,b, \rho, z_0\right)^{\frac{9-3 p}{6-p}} M^{\frac{3 p}{6-p}}.
\end{equation}
For $\text{case}\,2:3\leq p\leq 6$, $J_1$ can be showed as
\begin{equation}\label{2.15}
	\begin{split}
		J_{1} &\leq C\left\|u-u_{x,\rho}\right\|_{L^{p}}^{\frac{p}{p-2}}\left\|u-u_{x_0, \rho}\right\|_{L^{2}}^{\frac{2(p-3)}{p-2}}+ C\left\|b-b_{x_0,\rho}\right\|_{L^{p}}^{\frac{p}{p-2}}\left\|b-b_{x_0, \rho}\right\|_{L^{2}}^{\frac{2(p-3)}{p-2}}\\
		&\leq C\left(\left\|u-u_{x_0, \rho}\right\|_{L^{2}}^{\frac{2(p-3)}{p-2}}+\left\|b-b_{x_0, \rho}\right\|_{L^{2}}^{\frac{2(p-3)}{p-2}}\right)\left(\left\|u-u_{x_0,\rho}\right\|_{L^{p}}^{\frac{p}{p-2}}+\left\|b-b_{x_0,\rho}\right\|_{L^{p}}^{\frac{p}{p-2}}\right)\\
		&\leq C \rho A\left( \rho, z\right)^{\frac{p-3}{p-2}}\left(\frac{1}{\rho}\right)^{\frac{1}{p-2}}\left(\left\|u-u_{x_0,\rho}\right\|_{L^{p}}^{\frac{p}{p-2}}+\left\|b-b_{x_0,\rho}\right\|_{L^{p}}^{\frac{p}{p-2}}\right).
	\end{split}
\end{equation}
Together with estimates  $J_{1}$ and $J_{2}$ and integrating with respect to time from $t_0-r^{2}$ to $t_0$, we get
\begin{equation}\label{2.24}
	\begin{aligned}
		& \int_{t_0-r^{2}}^{t_0} \int_{B_{r}\left(x_0\right)}|u|^{3}+|b|^{3} d x d s \\
		\leq & C\left(\frac{r}{\rho}\right)^{3} \int_{t_0-\rho^{2}}^{t_0} \int_{B_{\rho}\left(x_0\right)}|u|^{3}+|b|^{3}d x d s+C \rho A\left( u,b,\rho, z_0\right)^{\frac{p-3}{p-2}} \int_{t_0-r^{2}}^{t_0} f_p^{\frac{p}{ p-2}}+g_p^{\frac{p}{p-2}} d s,
	\end{aligned}
\end{equation}
under the assumption \eqref{2.1}, we infer that
\begin{equation}\label{2.17}
	\begin{aligned}
		\int_{t_0-r^{2}}^{t_0} f_p^{\frac{p}{ p-2}}+g_p^{\frac{p}{p-2}} d s \leq& \int_{t-r^{2}}^{t} \left(f_p+g_p\right)^{\frac{p}{p-2}} d s\\
		=& \frac{p}{p-2} \int_{0}^{\infty} \sigma^{\frac{2}{p-2}}\left|\left\{s \in\left[t-r^{2}, t\right]: (f_p+g_p)(s)>\sigma\right\}\right| d \sigma \\
		=& \frac{p}{p-2}\left\{\int_{0}^{R} \sigma^{\frac{2}{p-2}}\left|\left\{s \in\left[t-r^{2}, t\right]: (f_p+g_p)(s)>\sigma\right\}\right| d \sigma\right.\\
		&+\int_{R}^{\infty} \sigma^{\frac{2}{p-2}}\left|\left\{s \in\left[t-r^{2}, t\right]: (f_p+g_p)(s)>\sigma\right\}\right| d \sigma \\
		\leq & \frac{p}{p-2} \int_{0}^{R} \sigma^{\frac{2}{p-2}} r^{2} d \sigma+\frac{p}{p-2} \int_{R}^{\infty} \sigma^{\frac{2}{p-2}-\frac{2 p}{p-2}} d \sigma M^{\frac{2 p}{p-2}} \\
		=& r^{2} R^{\frac{p}{p-2}}+\frac{p^{2}}{(p-2)^{2}} R^{\frac{-p}{p-2}} M^{\frac{2 p}{p-2}} \\
		=&\left(1+\frac{p^{2}}{(p-2)^{2}}\right) r M^{\frac{p}{p-2}},
	\end{aligned}
\end{equation}
where we have take $R=r^{-\frac{p-2}{p}} M$. Combining \eqref{2.16} and \eqref{2.17} it follows that
\begin{equation}\label{2.18}
	\begin{aligned}
		& \int_{t_0-r^{2}}^{t_0} \int_{B_{r}\left(x_0\right)}|u|^{3}+|b|^{3} d x d s \\
		\leq & C\left(\frac{r}{\rho}\right)^{3} \int_{t_0-\rho^{2}}^{t_0} \int_{B_{\rho}\left(x_0\right)}|u|^{3}+|b|^{3}d x d s+C\left(1+\frac{p^{2}}{(p-2)^{2}}\right)  \rho rA\left(  u,b,\rho, z_0\right)^{\frac{p-3}{p-2}}  M^{\frac{p}{p-2}}.
	\end{aligned}
\end{equation}
Therefore, we arrive at
\begin{equation}\label{2.27}
	C\left(u,b,r, z_0\right) \leq C \frac{r}{\rho} C\left(u,b,\rho, z_0\right)+C\left(\frac{\rho}{r}\right) A\left( u,b,\rho, z_0\right)^{\frac{p-3}{p-2}} M^{\frac{p}{p-2}}.
\end{equation}
For $\text{case}\,3:6< p\leq \infty$, $J_1$ can be showed as
\begin{equation}\label{2.20}
	\begin{split}
		J_{1} &\leq C\left\|u-u_{x_0,\rho}\right\|_{L^{p}}^{\frac{3}{2}}\left\|u-u_{x_0, \rho}\right\|_{L^{2}}^{\frac{3}{2}}\rho^{\frac{3(p-6)}{4 p}}+ C\left\|b-b_{x_0,\rho}\right\|_{L^{p}}^{\frac{3}{2}}\left\|b-b_{x_0, \rho}\right\|_{L^{2}}^{\frac{3}{2}}\rho^{\frac{3(p-6)}{4 p}}\\
		&\leq C\left(\left\|u-u_{x_0, \rho}\right\|_{L^{2}}^{\frac{3}{2}}+\left\|b-b_{x_0, \rho}\right\|_{L^{2}}^{\frac{3}{2}}\right)\left(\left\|u-u_{x,\rho}\right\|_{L^{p}}^{\frac{3}{2}}+\left\|b-b_{x_0,\rho}\right\|_{L^{p}}^{\frac{3}{2}}\right)\rho^{\frac{3(p-6)}{4 p}}\\
		&\leq C \rho A^{\frac{3}{4}}\left(u,b, \rho, z_0\right)\rho^{\frac{3(p-6)}{4 p}+\frac{3}{4}}\left(\left\|u-u_{x,\rho}\right\|_{L^{p}}^{\frac{3}{2}}+\left\|b-b_{x_0,\rho}\right\|_{L^{p}}^{\frac{3}{2}}\right)\\
		&\leq C \rho^{\frac{3}{4}+\frac{3}{2 p}+\frac{3(p-6)}{4 p}} A^{\frac{3}{4}}\left(u,b, \rho, z_0\right) \left(f_p+g_p\right)^{\frac{3}{2}}.
	\end{split}
\end{equation}
Together with estimates  $J_{1}$ and $J_{2}$ and integrating with respect to time from $t_0-r^{2}$ to $t_0$, we get
\begin{equation}\label{2.21}
	\begin{aligned}
		& \int_{t_0-r^{2}}^{t_0} \int_{B_{r}\left(x_0\right)}|u|^{3}+|b|^{3} d x d s \\
		\leq & C\left(\frac{r}{\rho}\right)^{3} \int_{t_0-\rho^{2}}^{t_0} \int_{B_{\rho}\left(x_0\right)}|u|^{3}+|b|^{3}d x d s+C \rho^{\frac{3p-6}{2p}} A\left( u,b,\rho, z_0\right)^{\frac{3}{4}}\int_{t_0-r^{2}}^{t_0} \left(f_p+g_p\right)^{\frac{3}{2}}(s) ds,
	\end{aligned}
\end{equation}
under the assumption \eqref{2.1}, it follows that
\begin{equation}\label{2.22}
	\begin{aligned}
		& \int_{t_0-r^{2}}^{t_0} \left(f_p+g_p\right)^{\frac{3}{2}}(s) ds\\
		=& \frac{3}{2} \int_{0}^{\infty} \sigma^{\frac{1}{2}}\left|\left\{s \in\left[t-r^{2}, t\right]: \left(f_p+g_p\right)^{\frac{3}{2}}(s)>\sigma\right\}\right| d \sigma \\
		=& \frac{3}{2}\left\{\int_{0}^{R} \sigma^{\frac{1}{2}}\left|\left\{s \in\left[t-r^{2}, t\right]:\left(f_p+g_p\right)^{\frac{3}{2}}(s)>\sigma\right\}\right| d \sigma\right.\\
		&+\int_{R}^{\infty} \sigma^{\frac{1}{2}}\left|\left\{s \in\left[t-r^{2}, t\right]: \left(f_p+g_p\right)^{\frac{3}{2}}(s)>\sigma\right\}\right| d \sigma \\
		\leq & \frac{3}{2} \int_{0}^{R} \sigma^{\frac{1}{2}} r^{2} d \sigma+\frac{3}{2} \int_{R}^{\infty} \sigma^{\frac{1}{2}-\frac{2 p}{p-2}} d \sigma M^{\frac{2 p}{p-2}} \\
		=& r^{2} R^{\frac{3}{2}}-\frac{3(p-2)}{p-6} R^{\frac{3}{2}-\frac{2p}{p-2}} M^{\frac{2 p}{p-2}} \\
		=&\left(1+\frac{3(p-2)}{6+p}\right) r^{2-\frac{3 p-6}{2 p}} M^{\frac{3}{2}},
	\end{aligned}
\end{equation}
where we have take $R=r^{\frac{2-p}{p}} M$. Substituting \eqref{2.22} into \eqref{2.21}, one has
\begin{equation}\label{2.31}
	\begin{aligned}
		& \int_{t_0-r^{2}}^{t_0} \int_{B_{r}\left(x_0\right)}|u|^{3}+|\nabla Q|^{3} d x d s \\
		\leq & C\left(\frac{r}{\rho}\right)^{3} \int_{t_0-\rho^{2}}^{t_0} \int_{B_{\rho}\left(x_0\right)}|u|^{3}+|b|^{3}d x d s+C \rho^{\frac{3p-6}{2p}}r^{2-\frac{3 p-6}{2 p}} A\left(u,b, \rho, z_0\right)^{\frac{3}{4}} M^{\frac{3}{2}},
	\end{aligned}
\end{equation}
which allows us to get
\begin{equation}\label{2.24}
	C\left(u,b, r,z_0\right) \leq C \frac{r}{\rho} C\left( u,b,\rho, z_0\right)+C\left(\frac{\rho}{r}\right)^{\frac{3}{2}} A\left(u,b,\rho, z_0\right)^{\frac{3}{4}} M^{\frac{3}{2}}.
\end{equation}
We complete the proof of Lemma \ref{le2.1} under the assumption \eqref{2.1}.

Now, we consider the assumption \eqref{2.2} is true, we denote $$f_{p}(t)=\sup _{\rho \leq 1}\left(\frac{1}{\rho} \int_{B_{\rho}\left(x_{0}\right)}|u(t, x)|^{p} d x\right)^{\frac{1}{p}}$$
and
$$g_{p}(t)=\sup _{\rho \leq 1}\left(\frac{1}{\rho} \int_{B_{\rho}\left(x_{0}\right)}|b(t, x)|^{p} d x\right)^{\frac{1}{p}},$$
and we can modify the processes of proofs as follows:

In the case $2 \leq p<3$, we can replace $\left\|u-u_{\rho}\right\|_{L^{p}\left(B_{\rho}\right)}$ and $\left\|b-b_{\rho}\right\|_{L^{p}\left(B_{\rho}\right)}$ by $C_{p}\|u\|_{L^{p}\left(B_{\rho}\right)}$ and
$C_{p}\|b\|_{L^{p}\left(B_{\rho}\right)}$ in \eqref{2.8} and repeat the processes of proofs for \eqref{2.9}-\eqref{2.13} to get \eqref{2.3}. The difference is that in this case, the constant $C_{p}$ is depended on $p$. Noticing $2 \leq p<3$, we can choose a large enough constant $C$ to get rid of the dependence on $p$.

If $3 \leq p \leq 6$ or $6<p \leq \infty$, we just need to replace $\int_{B_{r}\left(x_{0}\right)}(\left|u-u_{x_{0}, \rho}\right|^{3}+ \left|b-b_{x_{0}, \rho}\right|^{3})d x$ by $C \int_{B_{\rho}\left(x_{0}\right)}(|u|^{3}+|b|^{3}) d x$ and repeat the processes of proofs step by step to get \eqref{2.4} and \eqref{2.5}. Thus the proof of Lemma \ref{le2.1} is completed under the assumption \eqref{2.2}.

\begin{lemma}\label{le2.2}
	Let $(u, b, P)$ be a suitable weak solution of equations \eqref{1.1} on $Q_1(z_0)$ satisfying
\begin{equation}
	\begin{aligned}\label{2.25}
		&\left\| \left(\sup _{\rho\leq 1}\left(\frac{1}{\rho} \int_{B_{\eta}\left(x_{0}\right)}\left|u(x, t)-u_{x_{0},\rho}\right|^{p} d x\right)^{\frac{1}{p}},  \sup _{\rho \leq 1}\left(\frac{1}{\rho} \int_{B_{\eta}\left(x_{0}\right)}\left|b(x, t)-b_{x_{0}, \rho}\right|^{p} d x\right)^{\frac{1}{p}}\right) \right\|_{L^{q, \infty}\left[t_{0}-1, t_{0}\right]}\\
		& =M < \infty	
	\end{aligned}
\end{equation}
or
\begin{equation}
	\begin{aligned}\label{2.26}
		\left\|\left(\sup _{\rho \leq 1}\left(\frac{1}{\rho} \int_{B_{\rho}\left(x_{0}\right)}|u(x, t)|^{p} d x\right)^{\frac{1}{p}}, \sup _{\eta \leq 1}\left(\frac{1}{\rho} \int_{B_{\rho}\left(x_{0}\right)}|b(x, t)|^{p} d x\right)^{\frac{1}{p}}\right)\right\|_{L^{q, \infty}\left[t_{0}-1, t_{0}\right]}  =M < \infty,
	\end{aligned}
\end{equation}
where $p$ and $q$ satisfy $\frac{1}{q}+\frac{1}{p}=\frac{1}{2}$,	then there exists a constant $\rho_{0}>0$ only depended on $A\left(u, b, 1, z_{0}\right)$, $B\left(u, b, 1, z_{0}\right)$, $C\left(u, b, 1, z_{0}\right)$ and $D\left(P, 1, z_{0}\right)$, such that for $r \leq \rho_{0}$, it follows
$$
A(u,b, r, z_0)+B(u,b, r, z_0)+C(u,b, r, z_0)+D(P, r, z_0) \leq C(M).
$$
\end{lemma}
\textbf{Proof.}\quad Without loss of generality, we set $z_0=0$ and for any $2r<\rho<1$, let $\zeta$ be a cutoff function, which vanishes outside of $Q_{\rho}$ and equals 1 in $Q_{\frac{\rho}{2}}$, and satisfies
$$
|\nabla\zeta|\leq C\rho^{-1},\quad |\partial_t\zeta|, |\Delta\zeta|\leq C\rho^{-2}.
$$
Define the backward heat kernel as
$$
\Gamma(x,t)=\frac{1}{4\pi(r^2-t)^{\frac{3}{2}}}e^{-\frac{|x|^2}{4(r^2-t)}}.
$$
Taking the test function $\phi=\Gamma\zeta$ in the local energy inequality, and noting that
$$
(\partial_t+\Delta)\Gamma=0,
$$
we have
\begin{equation}
	\begin{split}\label{2.14}
		&\sup\limits_t\int_{B_r} (|u(x, t)|^2 +|b(x, t)|^2)\phi d x+\int_{Q_r}(|\nabla u(x, t)|^2 +|\nabla b(x, t)|^2)\phi dx dt\\
		\leq &\int_{Q_\rho} \left[(|u|^2 +|b|^2)(\Delta\phi+\partial_t\phi)+u\cdot\nabla\phi(|u|^2 +|b|^2+2p)-
		(b\cdot u)(b\cdot\nabla\phi)\right]dxdt\\
		\leq &\int_{Q_\rho} \left[(|u|^2 +|b|^2)(\Gamma\Delta\zeta+\Gamma\partial_t\zeta+2\nabla\Gamma\cdot\nabla\zeta)+|\nabla\phi||u|(|u|^2 +|b|^2+2P)\right]dxdt.
	\end{split}
\end{equation}
Direct computation shows that
$$
Cr^{-3}\geq\Gamma(x,t)\geq C^{-1}r^{-3}\quad \text{in}\quad Q_r;
$$
$$
|\nabla\phi|\leq |\nabla\Gamma|\zeta+\Gamma|\nabla\zeta|\leq Cr^{-4};
$$
$$
|\Gamma\Delta\zeta|+|\Gamma\partial_t\zeta|+2|\nabla\Gamma\cdot\nabla\zeta|\leq C\rho^{-5}.
$$
From these inequalities and the H\"older's inequality, it follows that
\begin{equation}
	\begin{split}\label{2.28}
		&A(u,b,r)+B(u,b,r)\\
		\leq &r^2\int_{Q_\rho} \left[(|u|^2 +|b|^2)(\Gamma\Delta\zeta+\Gamma\partial_t\zeta+2\nabla\Gamma\cdot\nabla\zeta)+|\nabla\phi||u|(|u|^2 +|b|^2+2P)\right]dxdt\\
		\leq &C\left(\frac{r}{\rho}\right)^2A\left(u,b,\rho\right)+C\left(\frac{\rho}{r}\right)^2\rho^{-2}\int_{Q_\rho} |u|^3 +|b|^3+2|u||p|dxdt\\
		=&C\left(\frac{r}{\rho}\right)^2A\left(u,b,\rho\right)+C\left(\frac{\rho}{r}\right)^2C\left(u,b,\rho\right)
		+C\left(\frac{\rho}{r}\right)^2\rho^{-2}\int_{Q_\rho} 2|u||P|dxdt\\
		\leq&C\left(\frac{r}{\rho}\right)^2A\left(u,b,\rho\right)+C\left(\frac{\rho}{r}\right)^2C\left(u,b,\rho\right)
		+C\left(\frac{\rho}{r}\right)^2C^{\frac{1}{3}}\left(u, b,\rho\right)D^{\frac{2}{3}}(P, \rho)\\
		\leq&C\left(\frac{r}{\rho}\right)^2A\left(u,b,\rho\right)+C\left(\frac{\rho}{r}\right)^2C\left(u,b,\rho\right)
		+C\left(\frac{\rho}{r}\right)^2D(P,\rho).
	\end{split}
\end{equation}
Now,we turn to show bounds on $D(P, r)$. Note that $p$ satisfies the following equation in distribution sense:
$$
-\Delta P=\partial_i\partial_j\left(u_iu_j-b_ib_j\right),$$
which in turn can be derived by calculating the divergence of $\eqref{1.1}_1$. Let $\eta(x)$ be a cut-off function, which equals 1 in $B_{\frac{\rho}{2}}$ and
vanishes outside of $B_\rho$. Set $P=P_1+P_2$ with
$$
P_1=\frac{1}{4\pi}\int_{\mathbb{R}^3}\frac{1}{|x-y|}\partial_i\partial_j[\left(u_iu_j-b_ib_j\right)\eta]dy
$$
and $P_2$ is harmonic in $B_{\frac{\rho}{2}}$. Using the Calder\'on-Zygmund inequality, one know that
$$\int_{B_{\rho}}|P_1|^{\frac{3}{2}} dx\leq C\int_{B_{\rho}}|u|^3+|b|^3dx.$$
Since $p_2$ is harmonic in $B_{\frac{\rho}{2}}$, we conclude that
$$\int_{B_{r}}|P_2|^{\frac{3}{2}} dx\leq Cr^3\sup\limits_{x\in B_{r}}|P_2(x)|^{\frac{3}{2}}\leq C\left(\frac{r}{\rho}\right)^3\int_{B_{\rho}}|P_2(x)|^{\frac{3}{2}}dx$$
for any $r\leq\frac{\rho}{2}$.
Hence, we infer that
\begin{equation*}
	\begin{split}
		\int_{B_{r}}|P|^{\frac{3}{2}} dx&\leq C\int_{B_{\rho}}|u|^3+|b|^3dx+C\left(\frac{r}{\rho}\right)^3\int_{B_{\rho}}|P-P_1|^{\frac{3}{2}}dx\\
		&\leq C\int_{B_{\rho}}|u|^3+|b|^3dx+C\left(\frac{r}{\rho}\right)^3\int_{B_{\rho}}|P|^{\frac{3}{2}}dx.
	\end{split}
\end{equation*}
Integrating above inequality with respect to $t$ from $-r^2$ to $0$, we get
\begin{equation*}
	\begin{split}
		\int_{Q_{r}}|P|^{\frac{3}{2}} dxdt
		&\leq C\int_{Q_{\rho}}|u|^3+|b|^3dxdt+C\left(\frac{r}{\rho}\right)^3\int_{Q_{\rho}}|P|^{\frac{3}{2}}dxdt,
	\end{split}
\end{equation*}
which implies
\begin{equation}
	\begin{split}\label{2.29}
		D(P,r)
		&\leq C\frac{1}{r^2}\int_{Q_{\rho}}\left(|u|^3+|b|^3\right)dxdt+C\frac{1}{r^2}\left(\frac{r}{\rho}\right)^3\int_{Q_{\rho}}|P|^{\frac{3}{2}}dxdt\\
		&\leq C\left(\frac{\rho}{r}\right)^2C(u,b,\rho)+C\frac{r}{\rho}D(P,\rho).
	\end{split}
\end{equation}
By using Young's inequality, we calculate
\begin{equation}
	\begin{split}\label{2.30}
		D^{\frac{8}{7}}(P,r)
		&\leq C\left(\frac{\rho}{r}\right)^{\frac{16}{7}}C^{\frac{8}{7}}(u,b,\rho)+C\left(\frac{r}{\rho}\right)^{\frac{8}{7}}D^{\frac{8}{7}}(P,\rho)\\
		&\leq C\left(\frac{\rho}{r}\right)^{\frac{16}{7}}\left(\frac{\rho}{r}\right)^{\frac{8}{7}}\left(\frac{r}{\rho}\right)^{\frac{8}{7}}C^{\frac{8}{7}}(u,b,\rho)
		+C\left(\frac{r}{\rho}\right)^{\frac{8}{7}}D^{\frac{8}{7}}(P,\rho)\\
		&\leq \left(\frac{r}{\rho}\right)^{\frac{7}{6}}C^{\frac{7}{6}}(u,b,\rho)+C\left(\frac{r}{\rho}\right)^{\frac{8}{7}}D^{\frac{8}{7}}(P,\rho)+C\left(\frac{\rho}{r}\right)^{168}.
	\end{split}
\end{equation}
Similarly, we show some crucial bounds for $C(u,b,r)$ under the assumptions of Lemma \ref{le2.1}. First, from Lemma \ref{le2.1}, we know that
\begin{equation*}\left\{\begin{array}{ll}
		C\left( u,b,r\right) \leq C \frac{r}{\rho} C\left(u,b,\rho\right)+C\left(\frac{\rho}{r}\right)^{2} B\left(u,b,\rho\right)^{\frac{9-3 p}{6-p}} M^{\frac{3 p}{6-p}}, \,\,2 \leq p<3;\\
		C\left(u,b,r\right) \leq C \frac{r}{\rho} C\left(u,b,\rho\right)+C\left(\frac{\rho}{r}\right) A\left(u,b,\rho\right)^{\frac{p-3}{p-2}} M^{\frac{p}{p-2}},\,\, 3 \leq p \leq 6; \\
		C\left( u,b,r\right) \leq C \frac{r}{\rho} C\left(u,b, \rho\right)+C\left(\frac{\rho}{r}\right)^{\frac{3}{2}} A\left(u,b,\rho\right)^{\frac{3}{4}} M^{\frac{3}{2}}, \,\,6<p \leq \infty.
	\end{array}\right.\end{equation*}
Clearly, we have
\begin{equation*}\left\{\begin{array}{ll}
		C^{\frac{7}{6}}\left(u,b, r\right)  \leq C\left(\frac{r}{\rho}\right)^{\frac{7}{6}} C^{\frac{7}{6}}(u,b, \rho)+C\left(\frac{r}{\rho}\right)^{\frac{7}{6}} B(u,b, \rho)+C\left(\frac{\rho}{r}\right)^{40} M^{\frac{7 p}{5 p-9}}, \,\,2 \leq p<3;\\
		C^{\frac{7}{6}}\left(u,b, r\right)  \leq C\left(\frac{r}{\rho}\right)^{\frac{7}{6}} C^{\frac{7}{6}}(u,b, \rho)+\left(\frac{r}{\rho}\right) A(u,b, \rho)+C\left(\frac{\rho}{r}\right)^{4} M^{\frac{7 p}{9-p}},\,\, 3 \leq p \leq 6; \\
	C^{\frac{7}{6}}\left(u,b, r\right)  \leq C\left(\frac{r}{\rho}\right)^{\frac{7}{6}} C^{\frac{7}{6}}(u,b, \rho)+\left(\frac{r}{\rho}\right)^{\frac{7}{6}} A(u,b, \rho)+C\left(\frac{\rho}{r}\right)^{23} M^{14}, \,\,6<p \leq \infty.
	\end{array}\right.\end{equation*}
From \eqref{2.28} and combine with above estimates of $C^{\frac{7}{6}}(r)$ and $D^{\frac{8}{7}}(r)$, we discover that
\begin{equation}
	\begin{split}\label{2.31}
		&A(r)+B(r)+C^{\frac{7}{6}}(r)+D^{\frac{8}{7}}(r)\\
		\leq&C\left(\frac{r}{\rho}\right)\left(A\left(\rho\right)+B\left(\rho\right)\right)+C\left(\frac{r}{\rho}\right)^{\frac{7}{6}}C^{\frac{7}{6}}(\rho)
		+C\left(\frac{r}{\rho}\right)^{\frac{8}{7}}D^{\frac{8}{7}}(\rho)+C\left(\frac{\rho}{r}\right)^{40} M^{\frac{7 p}{5 p-9}}\\
		&+C\left(\left(\frac{\rho}{r}\right)^{14}+\left(\frac{\rho}{r}\right)^{16}+\left(\frac{\rho}{r}\right)^{168}\right),
	\end{split}
\end{equation}
\begin{equation}
	\begin{split}\label{2.32}
		&A(r)+B(r)+C^{\frac{7}{6}}(r)+D^{\frac{8}{7}}(r)\\
		\leq&C\left(\frac{r}{\rho}\right)\left(A\left(\rho\right)+B\left(\rho\right)\right)+C\left(\frac{r}{\rho}\right)^{\frac{7}{6}}C^{\frac{7}{6}}(\rho)
		+C\left(\frac{r}{\rho}\right)^{\frac{8}{7}}D^{\frac{8}{7}}(\rho)+C\left(\frac{\rho}{r}\right)^{4} M^{\frac{7 p}{9-p}}\\
		&+C\left(\left(\frac{\rho}{r}\right)^{14}+\left(\frac{\rho}{r}\right)^{16}+\left(\frac{\rho}{r}\right)^{168}\right),
	\end{split}
\end{equation}
\begin{equation}
	\begin{split}\label{2.33}
		&A(r)+B(r)+C^{\frac{7}{6}}(r)+D^{\frac{8}{7}}(r)\\
		\leq&C\left(\frac{r}{\rho}\right)\left(A\left(\rho\right)+B\left(\rho\right)\right)+C\left(\frac{r}{\rho}\right)^{\frac{7}{6}}C^{\frac{7}{6}}(\rho)
		+C\left(\frac{r}{\rho}\right)^{\frac{8}{7}}D^{\frac{8}{7}}(\rho)+C\left(\frac{\rho}{r}\right)^{23} M^{14}\\
		&+C\left(\left(\frac{\rho}{r}\right)^{14}+\left(\frac{\rho}{r}\right)^{16}+\left(\frac{\rho}{r}\right)^{168}\right).
	\end{split}
\end{equation}
These estimates combined with the fact that $\frac{r}{\rho}<1$, one get
\begin{equation*}\left\{\begin{array}{ll}
		G(r)\leq C\left(\frac{r}{\rho}\right) G(\rho)+C\left(1+M^{\frac{7 p}{5 p-9}}\right)\left(\frac{\rho}{r}\right)^{168},\\
		G(r)\leq C\left(\frac{r}{\rho}\right) G(\rho)+C\left(1+M^{\frac{7 p}{9-p}}\right)\left(\frac{\rho}{r}\right)^{168},
		\\
		G(r)\leq C\left(\frac{r}{\rho}\right) G(\rho)+C\left(1+M^{14}\right)\left(\frac{\rho}{r}\right)^{168}.
	\end{array}\right.\end{equation*}
Let $r=\theta\rho$ with $\theta<\frac{1}{2}$ such that $C\cdot \frac{r}{\rho}<\frac{1}{2}$, i.e.,
\begin{equation}\label{2.34}
	\left\{\begin{array}{ll}
		G(\theta\rho)\leq \frac{1}{2} G(\rho)+C\left(1+M^{\frac{7 p}{5 p-9}}\right)\theta^{-168},\\
		G(\theta\rho)\leq \frac{1}{2} G(\rho)+C\left(1+M^{\frac{7 p}{9-p}}\right)\theta^{-168},
		\\
		G(\theta\rho)\leq \frac{1}{2} G(\rho)+C\left(1+M^{14}\right)\theta^{-168}.
	\end{array}\right.\end{equation}
Iterating the inequality \eqref{2.34} $k$ times yields
\begin{equation}\label{2.35}
	\left\{\begin{array}{ll}
		G(\theta^k\rho)\leq\frac{1}{2^k} G(\rho)+C\left(1+M^{\frac{7 p}{5 p-9}}\right)\theta^{-168},\\
		G(\theta^k\rho)\leq \frac{1}{2^k} G(\rho)+C\left(1+M^{\frac{7 p}{9-p}}\right)\theta^{-168},
		\\
		G(\theta^k\rho)\leq \frac{1}{2^k} G(\rho)+C\left(1+M^{14}\right)\theta^{-168}.
	\end{array}\right.\end{equation}
These complete the proof of Lemma \ref{le2.2}.
\begin{lemma}\label{le2.3}
	Let $(u, b, P)$ be a suitable weak solution of equations \eqref{1.1} on $Q_1(z_0)$. For any fixed $\epsilon>0$, there exists two constants $\delta$ and $r^*$ depended $\epsilon$ on  such that if
\begin{equation}
	\begin{aligned}\label{2.36}
		&\left\| \left(\sup _{\rho\leq 1}\left(\frac{1}{\rho} \int_{B_{\rho}\left(x_{0}\right)}\left|u(x, t)-u_{x_{0},\rho}\right|^{p} d x\right)^{\frac{1}{p}},  \sup _{\rho \leq 1}\left(\frac{1}{\rho} \int_{B_{\rho}\left(x_{0}\right)}\left|b(x, t)-b_{x_{0}, \rho}\right|^{p} d x\right)^{\frac{1}{p}}\right) \right\|_{L^{q, \infty}\left[t_{0}-1, t_{0}\right]}\\
		& \leq \delta\ll 1	
	\end{aligned}
\end{equation}
or
\begin{equation}
	\begin{aligned}\label{2.37}
		\left\|\left(\sup _{\rho \leq 1}\left(\frac{1}{\rho} \int_{B_{\rho}\left(x_{0}\right)}|u(x, t)|^{p} d x\right)^{\frac{1}{p}}, \sup _{\eta \leq 1}\left(\frac{1}{\rho} \int_{B_{\rho}\left(x_{0}\right)}|b(x, t)|^{p} d x\right)^{\frac{1}{p}}\right)\right\|_{L^{q, \infty}\left[t_{0}-1, t_{0}\right]}  \leq \delta\ll 1,
	\end{aligned}
\end{equation}
where $p$ and $q$ satisfy $\frac{1}{q}+\frac{1}{p}=\frac{1}{2}$,	then it follows
$$
C\left(u, b, r, z_{0}\right)+D\left(P, r, z_{0}\right) \leq \epsilon,\quad 0<r \leq r^{*}<1.
$$
\end{lemma}
\textbf{Proof.} Without loss of generality, we assume $z_0 = 0$ and $\delta \ll 1$. From Lemma \ref{le2.2}, we know that if $\rho \leq \rho_{0}$, then
\begin{equation}\label{2.38}
A( u,b,\rho)+B(u,b, \rho)+C( u,b,\rho)+ D( P,\rho)\leq C.
\end{equation}
Therefore, we choosing $\rho \leq \rho_{0}$ in \eqref{2.3} $(2\leq p<3)$ and using \eqref{2.38} to obtain
\begin{equation}
\begin{aligned}\label{2.39}
 C(u,b,r)
\leq & C \frac{r}{\rho} C\left(u,b, \rho\right)+C\left(\frac{\rho}{r}\right)^{2} B^{\frac{9-3 p}{6-p}}\left(u,b, \rho\right) \delta^{\frac{3 p}{6-p}} \\
\leq & C \frac{r}{\rho} +C \left(\frac{\rho}{r}\right)^{2}\delta^{\frac{3}{2}}.
\end{aligned}
\end{equation}
We take $r=2^{-k}\rho_0$ in \eqref{2.39} and obtain by a standard iterative argument that
\begin{equation}\label{2.40}
\begin{aligned}
	C\left(u, b, 2^{-k} \rho_{0}\right) & \leq \frac{1}{2} C\left(u, b, 2^{-k+1} \rho_{0}\right)+C 2^{2} \delta^{\frac{3}{2}} \\
	& \leq\left(\frac{1}{2}\right)^{k} C\left(u, b, \rho_{0}\right)+\sum_{j=0}^{k-1}\left(\frac{1}{2}\right)^{j} C\left(2^{2} \delta^{\frac{3}{2}}\right) \\
	& \leq\left(\frac{1}{2}\right)^{k} C\left(u, b, \rho_{0}\right)+C 2^{2} \delta^{\frac{3}{2}}.
\end{aligned}
\end{equation}
Similarly, using the same trick to see that
\begin{equation}\label{2.41}
	\begin{aligned}
		D\left(P, 2^{-k_{0} k} \rho_{0}\right) & \leq \frac{1}{2} D\left(P, 2^{-k_{0}(k-1)} \rho_{0}\right)+C 2^{2 k_{0}}\left(\left(\frac{1}{2}\right)^{k_{0}(k-1)} C\left(u, b, \rho_{0}\right)+C 2^{2} \delta^{\frac{3}{2}}\right) \\
		& \leq\left(\frac{1}{2}\right)^{k} D\left(P, 2^{-k_{0}} \rho_{0}\right)+C 2^{2 k_{0}}\left(\left(\frac{1}{2}\right)^{k_{0}(k-1)} C\left(u, b, \rho_{0}\right)+C 2^{2} \delta^{\frac{3}{2}}\right),
	\end{aligned}
\end{equation}
where we choose $k_{0}$ such that $C 2^{-k_{0}} \leq \frac{1}{2}$.
Summing up \eqref{2.40} and \eqref{2.41}, we have
\begin{equation}\label{2.42}
	\begin{aligned}
		& C\left(u, b, 2^{-k_{0} k} \rho_{0}\right)+D\left(P, 2^{-k_{0} k} \rho_{0}\right) \\
		\leq &\left(\left(\frac{1}{2}\right)^{k}+\left(\frac{1}{2}\right)^{k_{0} k}\right)\left(C\left(u, b, \rho_{0}\right)+D\left(P, 2^{-k_{0}} \rho_{0}\right)\right)+C 2^{2 k_{0}}\left(\left(\frac{1}{2}\right)^{k_{0}(k-1)} C\left(u, b, \rho_{0}\right)+C 2^{2} \delta^{\frac{3}{2}}\right) \\
		\leq & C\left(\left(\frac{1}{2}\right)^{k}+\left(\frac{1}{2}\right)^{k_{0} k}\right)+C 2^{2 k_{0}}\left(\left(\frac{1}{2}\right)^{k_{0}(k-1)}+C 2^{2} \delta^{\frac{3}{2}}\right),
	\end{aligned}
\end{equation}
and then we take $k$ large enough such that $C\left(2^{-k}+2^{-k k_{0}}+2^{-k k_{0}+3 k_{0}}\right) \leq \frac{\epsilon}{2}$ and $r^{*}=2^{-k k_{0}} \rho_{0}$, then take $\delta \leq\left(\frac{\epsilon}{2^{2 k_{0}+3} C}\right)^{\frac{2}{3}}$, we discover that
$$
C\left(u, b, r\right)+D\left(P, r\right) \leq \epsilon,\quad 0<r \leq r^{*}<1.
$$
Next, by choosing $\rho\leq\rho_{0}$ in \eqref{2.4} $(3 \leq p \leq 6)$ and using \eqref{2.38}, we infer that
	\begin{equation}
		\begin{aligned}\label{2.43}
	C\left(u,b, r\right) \leq &C \frac{r}{\rho} C\left(u,b, \rho_0\right)+C\left(\frac{\rho}{r}\right) A^{\frac{p-3}{p-2}}\left(u,b, \rho\right) \delta^{\frac{p}{p-2}}\\
	\leq & C \frac{r}{\rho} +C \left(\frac{\rho}{r}\right)\delta^{\frac{3}{2}}\\
	\leq & C \frac{r}{\rho} +C \left(\frac{\rho}{r}\right)^2\delta^{\frac{3}{2}},
\end{aligned}
\end{equation}
where we have used $\frac{p-3}{p-2} \leq 1$ and $\frac{p}{p-2} \geq \frac{3}{2}$.
By choosing $\rho\leq\rho_{0}$ in \eqref{2.5} $(6<p \leq \infty)$ and using \eqref{2.38}, we obtain by similar computations that
	\begin{equation}\label{2.44}
C(u,b, r) \leq C \frac{r}{\rho}+C\left(\frac{\rho}{r}\right)^{2} \delta^{\frac{3}{2}}.
\end{equation}
Therefore, repeating the process of \eqref{2.40}-\eqref{2.42}, we also have
$$
C\left(u, b, r\right)+D\left(P, r\right) \leq \epsilon,\quad 0<r \leq r^{*}<1.
$$
These impliy that the Lemma \ref{le2.3} holds true.
\begin{lemma}\label{le2.4} \cite{TAN1} For any given $r>0$ and $2 \leq p \leq \infty$, it follows that
	\begin{equation}\label{2.45}
	\left(\sup _{\eta \leq r} \frac{1}{\eta} \int_{B_{\eta}}|f|^{p} d x\right)^{\frac{1}{p}} \leq C\|f\|_{L^{\frac{3 p}{2}, \infty}\left(B_{r}\right)},
\end{equation}
	where $C>0$ is a constant independent on $p$.
\end{lemma}
\begin{lemma}\cite{HX1}\label{le2.5}\,\, Let $(u, b, P)$ be a suitable weak solution of the MHD equations
	\eqref{1.1} in $Q_{1}\left(z_{0}\right)$. There exists an $\varepsilon>0$ such that if
	$$
	C(u, b, r,z_0)+D(P, r,z_0) \leq \varepsilon
	$$
	for some $r>0$, then $(u, b)$ is regular at $z_{0}$
	.\end{lemma}
\section{Proof of main results}
This section is devoted to the proof of the main theorems. We first prove the local in
space regularity near initial time, i.e. Theorem  \ref{th1}. Then, we prove the concentration result,
i.e. Theorem \ref{th2}.
\subsection{Proof of Theorem \ref{th1}}
From Lemma \ref{le2.3}, there exists two constants $\delta$ and $r^*$ depended $\epsilon$ on  such that if
\begin{equation}
	\begin{aligned}\label{3.1}
		&\left\| \left(\sup _{\eta\leq 1}\left(\frac{1}{\eta} \int_{B_{\eta}\left(x_{0}\right)}\left|u(x, t)-u_{x_{0},\rho}\right|^{p} d x\right)^{\frac{1}{p}},  \sup _{\eta \leq 1}\left(\frac{1}{\eta} \int_{B_{\eta}\left(x_{0}\right)}\left|b(x, t)-b_{x_{0}, \eta}\right|^{p} d x\right)^{\frac{1}{p}}\right) \right\|_{L^{q, \infty}\left[t_{0}-1, t_{0}\right]}\\
		& \leq \delta	
	\end{aligned}
\end{equation}
or
\begin{equation}
	\begin{aligned}\label{3.2}
		\left\|\left(\sup _{\eta \leq 1}\left(\frac{1}{\eta} \int_{B_{\eta}\left(x_{0}\right)}|u(x, t)|^{p} d x\right)^{\frac{1}{p}}, \sup _{\eta \leq 1}\left(\frac{1}{\eta} \int_{B_{\eta}\left(x_{0}\right)}|b(x, t)|^{p} d x\right)^{\frac{1}{p}}\right)\right\|_{L^{q, \infty}\left[t_{0}-1, t_{0}\right]}  \leq \delta,
	\end{aligned}
\end{equation}
where $p$ and $q$ satisfy $\frac{1}{q}+\frac{1}{p}=\frac{1}{2}$,	then it follows that
$$
C\left(u, b, r, z_{0}\right)+D\left(P, r, z_{0}\right)  \leq \epsilon,\quad 0<r\leq r^{*}<1.
$$
By Lemma \ref{le2.5}, we deduce that $z_0$ is a regular point. This completes the proof of Theorem \ref{th1}.
\subsection{Proof of Theorem \ref{th2}}
The proof of Theorem \ref{th2} is by contradiction.
If Theorem \ref{th2}  is false, then there exists some $r_0 \in(0,1)$ such that it holds that
\begin{equation}\label{3.3}
\limsup _{t \rightarrow t_{0}^{-}}\left(\left\|u(t, x)-u(t)_{x_{0}, r_0}\right\|_{L^{3, \infty}\left(B_{r_0}\left(x_{0}\right)\right)}+\left\|b(t, x)-b(t)_{x_{0}, r_0}\right\|_{L^{3, \infty}\left(B_{r_0}\left(x_{0}\right)\right)}\right)\leq \delta^{*}
\end{equation}
or
\begin{equation}
	\begin{aligned}\label{3.4}
		\limsup _{t \rightarrow t_{0}^{-}}\left(t_{0}-t\right)^{\frac{1}{\mu_0}} r^{\frac{2}{\nu_0}-\frac{3}{p_0}}_0\|(u,b)(t)\|_{L^{p_0, \infty}\left(B_{r_0}\left(x_{0}\right)\right)}\leq\delta^{*}, \frac{1}{\mu_0}+\frac{1}{\nu_0}=\frac{1}{2} , 2 \leq \nu_0 \leq \frac{2p_0}{3}, 3<p_0\leq\infty
	\end{aligned}
\end{equation}
or
\begin{equation}
	\begin{aligned}\label{3.5}
	\limsup\limits _{t \rightarrow t_{0}^{-}}\left(t_{0}-t\right)^{\frac{1}{\mu_0}} r^{\frac{2}{\nu_0}-\frac{3}{p_0}+1}_0\|(\nabla u,\nabla b)(t)\|_{L^{p_0}\left(B_{r_0}\left(x_{0}\right)\right)}\leq\delta^{*} \text { for } \frac{1}{\mu_0}+\frac{1}{\nu_0}=\frac{1}{2},\, \nu_0 \in\left\{\begin{array}{ll}
		{[2, \infty],} & p_0\geq 3 \\
		{[2, \frac{2p_0}{3-p_0}],} & \frac{3}{2}\leq p_0<3
	\end{array}\right.
	\end{aligned}
\end{equation}
If \eqref{3.3} holds, by using the fact $$\min _{c \in \mathbb{R}} \int_{B_{\eta}\left(x_{0}\right)}|u-c|^{2} d x=\int_{B_{\eta}\left(x_{0}\right)}\left|u-u_{x_{0}, \eta}\right|^{2} d x,\quad
\min _{c \in \mathbb{R}} \int_{B_{\eta}\left(x_{0}\right)}|b-c|^{2} d x=\int_{B_{\eta}\left(x_{0}\right)}\left|b-b_{x_{0}, \eta}\right|^{2} d x$$
and Lemma \ref{le2.4}, we immediately get
\begin{equation}
	\begin{aligned}
		& \limsup _{t \rightarrow t_{0}} \sup _{\eta<r_{0}}\left(\frac{1}{\eta} \int_{B_{\eta}\left(x_{0}\right)}\left|u-u_{x_{0}, \eta}\right|^{2} d x\right)^{\frac{1}{2}}+\limsup _{t \rightarrow t_{0}} \sup _{\eta<r_{0}}\left(\frac{1}{\eta} \int_{B_{\eta}\left(x_{0}\right)}\left|b-b_{x_{0}, \eta}\right|^{2} d x\right)^{\frac{1}{2}} \\
		\leq & C \limsup _{t \rightarrow t_{0}} \sup _{\eta<r_{0}}\left(\eta^{-1}\left(\int_{B_{\eta}\left(x_{0}\right)}\left|u-u_{x_{0}, r_0}\right|^{2} d x\right)^{\frac{1}{2}}\right)+C \limsup _{t \rightarrow t_{0}} \sup _{\eta<r_{0}}\left(\eta^{-1}\left(\int_{B_{\eta}\left(x_{0}\right)}\left|b-b_{x_{0}, r_0}\right|^{2} d x\right)^{\frac{1}{2}}\right) \\
		\leq& C \limsup _{t \rightarrow t_{0}}\left\|u(t)-u(t)_{x_{0}, r_0}\right\|_{L^{3, \infty}\left(B_{r_{0}}\left(x_{0}\right)\right)}+C \limsup _{t \rightarrow t_{0}}\left\|b(t)-b(t)_{x_{0}, r_0}\right\|_{L^{3, \infty}\left(B_{r_{0}}\left(x_{0}\right)\right)} \\
		\leq & C \delta^{*}=\delta.
	\end{aligned}
\end{equation}
This implies that
	\begin{equation}
	\begin{aligned}\label{3.7}
		\left\|\left(\sup _{\eta \leq r_0^*}\left(\frac{1}{\eta} \int_{B_{\eta}\left(x_{0}\right)}\left|u(x,t)-u_{x_{0}, \eta}\right|^{2} d x\right)^{\frac{1}{2}}, \sup _{\eta \leq r_0^*}\left(\frac{1}{\eta} \int_{B_{\eta}\left(x_{0}\right)}\left|b(x,t)-b_{x_{0}, \eta}\right|^{2} d x\right)^{\frac{1}{2}}\right)\right\|_{L^{ \infty}\left[t_{0}-(r_0^*)^2, t_{0}\right]} \leq \delta
	\end{aligned}
\end{equation}
for some $r^*_0\leq r_0$. Due to Theorem \ref{th1}, we know that $z_{0}$ is a regular point, this is a contradiction.\\
If \eqref{3.4} holds, by using H\"older's inequality and Lemma \ref{le2.4}, we find that
\begin{equation}
	\begin{aligned}
		& \limsup _{t \rightarrow t_{0}}\left(t_{0}-t\right)^{\frac{1}{\mu_{0}}} \left(\sup _{\eta<r_{0}}\left(\frac{1}{\eta} \int_{B_{\eta}\left(x_{0}\right)}|u|^{\nu_{0}} d x\right)^{\frac{1}{\nu_{0}}}+\sup _{\eta<r_{0}}\left(\frac{1}{\eta} \int_{B_{\eta}\left(x_{0}\right)}|b|^{\nu_{0}} d x\right)^{\frac{1}{\nu_{0}}}\right) \\
		\leq & \limsup _{t \rightarrow t_{0}}\left(t_{0}-t\right)^{\frac{1}{\mu_{0}}} \sup _{\eta<r_{0}}\left(\left(\eta^{2-\frac{3 \nu_{0}}{ p_{0}}}\left(\frac{1}{\eta} \int_{B_{\eta}\left(x_{0}\right)}|u|^{\frac{2 p_{0}}{3}} d x\right)^{\frac{3 \nu_{0}}{2 p_{0}}}\right)^{\frac{1}{\nu_{0}}}+\left(\eta^{2-\frac{3 \nu_{0}}{ p_{0}}}\left(\frac{1}{\eta} \int_{B_{\eta}\left(x_{0}\right)}|b|^{\frac{2 p_{0}}{3}} d x\right)^{\frac{3 \nu_{0}}{2 p_{0}}}\right)^{\frac{1}{\nu_{0}}}\right)  \\
		\leq & \limsup _{t \rightarrow t_{0}}\left(t_{0}-t\right)^{\frac{1}{\mu_{0}}} \sup _{\eta<r_{0}}\left(\eta^{\frac{2}{\nu_{0}}-\frac{3}{p_{0}}}\left(\frac{1}{\eta} \int_{B_{\eta}\left(x_{0}\right)}|u|^{\frac{2 p_{0}}{3}} d x\right)^{\frac{3}{2 p_{0}}}\right)+\sup _{\eta<r_{0}}\left(\eta^{\frac{2}{\nu_{0}}-\frac{3}{p_{0}}}\left(\frac{1}{\eta} \int_{B_{\eta}\left(x_{0}\right)}|b|^{\frac{2 p_{0}}{3}} d x\right)^{\frac{3}{2 p_{0}}}\right) \\
		\leq & \limsup _{t \rightarrow t_{0}}\left(t_{0}-t\right)^{\frac{1}{\mu_{0}}} \sup _{\eta<r_{0}}\left(\eta^{\frac{2}{\nu_{0}}-\frac{3}{p_{0}}}\|(u,b)(t)\|_{L^{p_{0}, \infty}\left(B_{\eta}\left(x_{0}\right)\right)}\right) \\
		\leq & C\limsup _{t \rightarrow t_{0}}\left(t_{0}-t\right)^{\frac{1}{\mu_{0}}}r_0^{\frac{2}{\nu_{0}}-\frac{3}{p_{0}}}\|(u,b)(t)\|_{L^{p_{0}, \infty}\left(B_{\eta}\left(x_{0}\right)\right)}\\
		\leq& C\delta^*=\delta,
	\end{aligned}
\end{equation}
where $\delta$ is the same constant in Theorem \ref{th1} and we choose $\delta^{*}=\frac{\delta}{C}$. This implies
	\begin{equation}
	\begin{aligned}\label{3.6}
		\left\|\left(\sup _{\eta \leq r^*_0}\left(\frac{1}{\eta} \int_{B_{\eta}\left(x_{0}\right)}|u(x, t)|^{\nu_0} d x\right)^{\frac{1}{\nu_0}}, \sup _{\eta \leq r^*_0}\left(\frac{1}{\eta} \int_{B_{\eta}\left(x_{0}\right)}|b(x, t)|^{\nu_0} d x\right)^{\frac{1}{\nu_0}}\right)\right\|_{L^{\mu_0, \infty}\left[t_{0}-(r_0^*)^2, t_{0}\right]} \leq \delta
	\end{aligned}
\end{equation}
for some $r^*_0\leq r_0$. Due to facts of Theorem \ref{th1}, we know that $z_{0}$ is a regular point, which is a contradiction.\\
If \eqref{3.5} is true, by using Poincar\'e's inequality, we have
	\begin{equation}
\begin{aligned}
	& \limsup _{t \rightarrow t_{0}}\left(t_{0}-t\right)^{\frac{1}{\mu_{0}}} \left(\sup _{\eta<r_{0}}\left(\frac{1}{\eta} \int_{B_{\eta}\left(x_{0}\right)}\left|u-u_{x_{0}, \eta}\right|^{\nu_{0}} d x\right)^{\frac{1}{\nu_{0}}}+\sup _{\eta<r_{0}}\left(\frac{1}{\eta} \int_{B_{\eta}\left(x_{0}\right)}\left|b-b_{x_{0}, \eta}\right|^{\nu_{0}} d x\right)^{\frac{1}{\nu_{0}}}\right) \\
	\leq & C \limsup _{t \rightarrow t_{0}}\left(t_{0}-t\right)^{\frac{1}{\mu_{0}}} \sup _{\eta<r_{0}}\left(\left(\eta^{2-\frac{3\nu_{0}}{p_0}+\nu_0}\left(\int_{B_{\eta}\left(x_{0}\right)}|\nabla u|^{p_0} d x\right)^{\frac{\nu_{0}}{p_0}}\right)^{\frac{1}{\nu_{0}}}+\left(\eta^{2-\frac{3\nu_{0}}{p_0}+\nu_0}\left(\int_{B_{\eta}\left(x_{0}\right)}|\nabla b|^{p_0} d x\right)^{\frac{\nu_{0}}{p_0}}\right)^{\frac{1}{\nu_{0}}}\right) \\
	\leq& C \limsup _{t \rightarrow t_{0}}\left(t_{0}-t\right)^{\frac{1}{\mu_{0}}} r_{0}^{\frac{2}{\nu_0}-\frac{3}{p_0}+1}\|(\nabla u,\nabla b)\|_{L^{p_0}\left(B_{r_{0}}\left(x_{0}\right)\right)} \\
	\leq& C  \delta^{*}=\delta,
\end{aligned}
\end{equation}
which implies
	\begin{equation}
	\begin{aligned}
		\left\| \left(\sup _{\eta \leq r_{0}^{*}}\left(\frac{1}{\eta} \int_{B_{\eta}\left(x_{0}\right)}\left|u(x, t)-u_{x_{0}, \eta}\right|^{\nu_0} d x\right)^{\frac{1}{\nu_0}},  \sup _{\eta \leq r_{0}^{*}}\left(\frac{1}{\eta} \int_{B_{\eta}\left(x_{0}\right)}\left|b(x, t)-b_{x_{0}, \eta}\right|^{\nu_0} d x\right)^{\frac{1}{\nu_0}}\right) \right\|_{L^{\mu_0, \infty}\left[t_{0}-\left(r_{0}^{*}\right)^{2}, t_{0}\right]} \leq \delta 	
	\end{aligned}
\end{equation}
for some $r_{0}^{*} \leq r_{0}$
and $\nu_0 \in\left\{\begin{array}{ll}
	{[2, \infty],} & p_0\geq 3 \\
	{[2, \frac{2p_0}{3-p_0}],} & \frac{3}{2}\leq p_0<3
\end{array}\right.$.
By using Theorem \ref{th1}, we deduce $z_{0}$ is a regular point, this is also a contradiction, hence concludes the proof of Theorem \ref{th2}.
\section*{Acknowledgments}

The authors are supported by the Construct Program of the Key Discipline in Hunan Province and NSFC Grant No. 11871209.

\end{document}